\documentclass[12pt, reqno]{amsart}
\usepackage{amsfonts}
\usepackage{graphicx}
\usepackage{amsmath, amscd, amssymb}
\usepackage[latin5]{inputenc}

\begin{document}

\begin{center}
\thispagestyle{empty} 
\markboth{\textbf{ }}
{\textbf{The New Existence and Uniqueness Results}}

\noindent {\Large \textbf{The New Existence and Uniqueness Results for Complex Nonlinear Fractional \\[0pt]
Differential Equation }}

\bigskip

\textbf{M. Şan$^{1,*}$ and K.N. Soltanov$^{2}$ }
\end{center}

\vspace{.03cm}

\begin{center}
\noindent \textit{$^{1}$Department of Mathematics, Faculty of Science, Çank{\i}r{\i} 
Karatekin University, Tr-18100, Çank{\i}r{\i}, Turkey\\[0pt]
e-mail address: mufitsan@karatekin.edu.tr } \\[0pt]
\noindent \textit{$^{2}$ Department of Mathematics, Faculty of Science,
Hacettepe University, Beytepe, Ankara, Turkey \\[0pt]
e-mail address: soltanov@hacettepe.edu.tr }
\end{center}

\medskip

\vspace{.4cm}

\noindent \textbf{Abstract.} In this article, we obtain existence and uniqueness results to some problems involving complex nonlinear fractional differential equations (FDEs) in the closed unit disc of $ 
\mathbb{C}.$ By help of these results, we prove that some IVPs for some fractional differential equations  with Caputo or Riemann-Liouville derivative admit at least one local (or unique) solution continuous on a closed interval $\left[0,R\right]$ and real analytic on $\left(0,R\right),$ where $0<R\leq 1.$ 
 
\medskip


\medskip


\bigskip

\noindent {\textbf{\large 1.Introduction}}   

\medskip\noindent In the last decade, there has been a great attention to fractional calculus and, in particular, fractional differential equations. Much as there are many researches on FDE involving Riemann-Liouville, Caputo, Grünwald-Letkinov, Caputo-Fabrizio derivatives, in the recent years some researchers have tried to develop new definitions for fractional derivative to model the physical phenomena better. In addition to this,  there have been numerous studies (for example, \cite{Campos}, \cite{Campos2}, \cite{Li}, \cite{Nish}, \cite{Ortigueira}-\cite{Owa}) to generalize the fractional derivatives defined on the real line mentioned above to those on the complex plane. 

In this work, we use a complex generalization of Riemann-Liouville fractional, which is given by Definition 2.1 and in \cite{Samko} and, we 
first consider the following problem with this definition   
$$D_{z}^{a}u(z)= f\big( z,u(z)\big)  \   \  (z\in\mathbb{U}) \eqno{(1.1)}$$
$$u(0)  =  b,\ \ \  \ \ \ \ \  \ \  \ \ \  \  \ \ \  \ \   $$
where $0<a<1,$  $\mathbb{U}$ is the open unit disc, $b\in\mathbb{C},$ and $f$ will be given in the sequel.   

The problem (1.1) with the fractional derivative defined by Owa in \cite{Owa} and with $b=0$ was first studied by Ibrahim and Darus \cite{Dar}. They claimed, in Theorem 4.1-4.2, that the problem (1.1) has an analytic solution when $f$ is analytic. But, if the problem they considered has an analytic solution $u(z)$, then being multivalent of  $D_{z}^{a}(u(z))$ on $\mathbb{U}$ contradicts with the analyticity of $f(z,u(z))$ on $\mathbb{U}.$

We, on the other hand, investigate the problem (1.1) provided that the function $f$ satisfies the following conditions: 

\smallskip

 (I) $f(z,t)$ is analytic on $\mathbb{D}\times \mathbb{C}$ and continuous on $\mathbb{D}^{*}\times \mathbb{C},$ and $z^{a}f(z,t)$ is
analytic on $\mathbb{U} \times \mathbb{C}$ and continuous on $\overline{\mathbb{U}}\times \mathbb{C},$ 
where $\mathbb{D}=\big\{ z\in\mathbb{U}: -\pi<\arg z\leq\pi\big\}$ and  $\mathbb{D}^{\ast }=\big\{ z\in\overline{\mathbb{U}}: -\pi<\arg z\leq\pi\big\},$ 

\smallskip  

 (II) $z^{a}f(z,b)\big|_{z=0}=b/\Gamma(1-a),$ 

\smallskip \noindent and we show the  existence and uniqueness of the solution analytic on the open unit disc $\mathbb{U}$ and continuous on the boundary of $\mathbb{U}.$ 
Imposing a condition in addition to (I) and using Schwarz Lemma, we prove the existence of the desired solution for the problem (1.1). For the uniqueness of the solution, we impose the Lipschitz-type condition on the function $f$ in addition to (I)-(II) and we again use Schwarz Lemma. Using this Lemma provides us to prove the uniqueness of the desired solution of the considered problem for the function $f$ which belongs to a class of analytic function larger than that obtained in the earlier studies (see Remark 3.9 (i)). \

\smallskip Secondly, in this study, we can establish the existence and uniqueness of the solution for the problem 
$$D_{z}^{a}(u(z)-u(0))= f\big(z,u(z)\big)  \   \  (z\in\mathbb{U})   \eqno{(1.2)}$$
$$u(0)  =  b,\ \ \  \ \ \ \            $$ 
where $f$ satisfies the similar conditions to those for the problem (1.1), 

Furthermore, by help of the results obtained for the problems given above,  we can reveal some existence and uniqueness results for real analytic and continuous solution of 
$$\mathcal{D}^{a}u(x)= f\big(x,u(x)\big)  \   \  \eqno{(1.3)}$$
$$u(0)  =  b,\ \ \  \ \ \ \  \ \    $$
where $\mathcal{D}^{a}$ is Riemann-Liouville or Caputo derivative defined in the real line. More precisely, for some appropriate $f(x,y),$  $u(z)$ will be a solution of the problem (1.1) (or the problem (1.2)) with the function $f(z,t)$ analytic continuation of $f(x,y)$ such that real part of $u(z)$ will be a solution of the problem (1.3) with Riemann-Liouville (or with Caputo derivative). There has been many approaches and techniques such as Tonelli's approach, monotone iterative technique etc. used by researchers (for example, \cite{Bale2}, \cite{Laks}, \cite{Zhang}) to establish the existence and uniqueness of continuous solution of the problem (1.3) with continuous right-hand side. Here we, on the other hand, use the complex variable approach for the same goal. 

It needs to note that  Kai Diethelm investigated the existence of real analytic solution of this problem. As he has shown in Theorem 2.1, for the existence of analytic solution, the condition $f\big(x,u(x)\big)=0$ for all $x\in [0,T]$ must be satisfied. This indicates that the occurrence of an real analytic solution on $(0,T)$ (and also continuous on $[0,T]$) to an equation of (1.3) with analytic right-hand side is a rare event. At the same time, if we consider the conditions we will pose for the problem (1.3) in the sequel, then we can say that the existence of analytic solution for the problem (1.3) under these conditions is not rare event. 

On the other hand, the geometric properties of the solutions to ordinary  differential equations in the complex plane have been discussed in many papers (for example, \cite{Saitoh},\cite{San2}). In
the present paper, we briefly mention about the geometric properties of the solutions to complex nonlinear fractional differential equation. We show that the univalence of $z^{a}f(z,t)$ does not imply the univalence of the solution of the problem (1.1).  

From reasons mentioned above we assert that our main results are new, and general than those obtained earlier in the complex plane. Moreover, using of Schwarz Lemma  provides us to obtain better results. Therefore, we think that the existence and uniqueness theorems  we obtained for the above problems may make an important contribution to the area of fractional differential equation.  
 
\medskip\noindent{\textbf{\large 2.Preliminaries}} 

\smallskip\noindent In this section, for the main results, we present several definitions and preliminary results. We begin with the definitions of the fractional integral and derivative which are applicable to not only analytic functions  but also integrable functions and given in \cite{Samko}, since we are interested in the continuous functions as well as analytic functions in the considered problems.  

\smallskip \noindent \textbf{Definition 2.1.} Let the function $u(z)$ defined on a certain domain of complex plane containing the points $0$ and $z.$ Then, the fractional integral and derivative of order $a$ $\big(0<a<1\big)$ of $u(z)$ are defined, respectively, by 
\begin{equation}
I_{z}^{a}u\left( z\right) :=\frac{1}{\Gamma \left( a\right) }%
\ \int_{0}^{z}\frac{u\left( \zeta \right) }{\left( z-\zeta \right)
^{1-a}}d\zeta,  \tag{2.1}
\end{equation}
and 
\begin{equation}
D_{z}^{a}u\left( z\right) =\frac{1}{\Gamma \left( 1-a\right) 
}\frac{d}{dz}\int_{0}^{z}\frac{u\left( \zeta \right) }{\left( z-\zeta
\right) ^{a}}d\zeta   \tag{2.2}
\end{equation}
where the integrations are along the straight line interval connecting points $0$ and $z$ as a rule, and with the principal value
$$(z-\zeta)^{{1-a}}=\left|z-\zeta\right|^{1-a}e^{i(1-a)\arg(z)}, \  \ \ \arg(z)\in (-\pi,\pi]. $$  

Let $u(x)=\Re\left\{u(z)\right\}$ on the real line. Then we get 
$$\mathcal{D}^{a}u(x)=D_{z}^{a}\left(\Re\left\{u(z)\right\}\right) \ \  \text{on} \ \ \mathbb{R} ,\eqno{(2.3)} $$
where $\mathcal{D}^{a}$ is the well-known Riemann Liouville derivative. 

\smallskip
Moreover, $I_{z}^{a}$ possess the semigroup property, i.e, if $u(z)$ be locally integrable (continuous) in a domain $G,$ then for almost all
(for all) $z\in G$ the following equality holds: 
\begin{equation*}
I_{z}^{a}I_{z}^{\beta } u(z)= I_{z}^{a+\beta }u(z)
\text{ \ \ \ } (a> 0,\beta> 0).
\end{equation*}

The further information for these definitions and related definitions can be found in \cite{Samko}.  

\smallskip We investigate the solution of the above problems in the space $\mathcal{B}_{R}$ with $0<R\leq 1$ (see \cite{Hofmann}) defined in the following:  
 
\smallskip \noindent\textbf{Definition 2.2.}  $\mathcal{B}_{R}$ is denoted the
space of functions which are analytic on $\mathbb{U_{R}}:=\big\{ z\in\mathbb{C}: \left|z\right|<R\big\}$ and continuous on the
boundary of $\mathbb{U_{R}}$. $\mathcal{B}_{R}$ is a Banach space when endowed with
the supremum norm.  The space $\mathcal{B}_{R}^{0}$ is represented by $\mathcal{B}_{R}^{0}=\{u\in\mathcal{B}_{R}:u(0)=0\}$ and we set $\mathcal{B}=\mathcal{B}_{R}$  for $R:=1.$ 

\smallskip As a consequences of Arzel\'a-Ascoli Theorem given in the complex plane and Schauder fixed point theorem (see \cite{Gamelin} and \cite{Zeidler}) one can give the following theorem:  

\smallskip \noindent \textbf{Theorem 2.3.} Let $M$ be a close bounded convex subset of a Banach space $X:=\big\{u:G\to\mathbb{C} \ \text{continuous}  : G\subset\mathbb{C} \ \text{compact}\big\}.$  If $T:M\rightarrow M$ is a continuous operator and $T(M)$ is a equicontinuous set on $G,$  
then $T$ has a fixed point in $M.$  

\smallskip For proving not only the existence but also uniqueness of the solution for the above problems we use  Banach fixed point theorem \cite{Zeidler}:  

\smallskip \noindent \textbf{Theorem 2.4.} If $(X,d)$ is a complete metric space and $T:X\rightarrow X$ is a contraction mapping, i.e  there is a $\beta$ ($0\leq\beta<1$) such that for all $x,y\in X$ $$d(Tx,Ty)\leq\beta d(x,y),$$ is satisfied, then $T$ has a unique fixed point. \\

In next sections, we show the existence and uniqueness of the desired solution for the considered problems by the help of Schwarz Lemma \cite{Dettman}:  

\smallskip \noindent\textbf{Lemma 2.5.} If $u(z)$ is analytic on
$\mathbb{U}_{R}:=\big\{ z\in\mathbb{C}: \left|z\right|<R\big\}$ and satisfies conditions $u(0)=0,$ $\left|u(z)\right|\leq r $ on $\mathbb{U}_{R}$, then 
\begin{align*}
\left|u(z)\right|\leq \frac{r}{R}\left|z\right|,   \  \  (\forall z\in\mathbb{U}_{R}) . \tag{2.4}
\end{align*}
\noindent In addition to above hypothesis, if $u$ is continuous
on $\overline{\mathbb{U}}_{R}$ and satisfies $\left|u(z)\right|\leq r $ for $\overline{\mathbb{U}}_{R},$
then the inequality (2.4) is satisfied for $\overline{\mathbb{U}}_{R}.$ \\
 
\medskip\noindent{\textbf{\large 3. Existence and Uniqueness Results Problems for (1.1) and (1.2)}} 

\smallskip\noindent We begin with proving the compositional relations in the following, which help us to define the equivalent form of the solution for the problem (1.1). 

\smallskip\noindent \textbf{Lemma 3.1.} Let $0<a<1$. Suppose that $u$ is a continuous and integrable function in $\mathbb{D}_{R}:=\big\{ z\in\mathbb{U}_{R}: -\pi<\arg z\leq\pi\big\}$ for an arbitrary fixed $R>0,$ then the fractional differential equation 
\begin{equation}
D_{z}^{a}u(z)=0 \  \  \  \left(z\in \mathbb{D}_{R}\right)  \tag{3.1}
\end{equation}
has the solutions which are only in the form $u(z)=cz^{a-1}$ with $c\in\mathbb{C}.$ 

\smallskip \noindent \textbf{Proof.}  For the proof, it is shown that there is a contradiction. It is obvious that $u(z)=cz^{a-1}$ ($c\in \mathbb{C}$) are the solutions of (3.1). We suppose that there exists a different solution $v$ of (3.1). Hence, $v(z)-cz^{a-1}$
are also the solutions of (3.1), since $D_{z}^{a}$ is a linear operator. By using $D_{z}^{a}=DI_{z}^{1-a}$ in (3.1), we have 
\begin{equation*}
I_{z}^{1-a}(v(z)-cz^{a-1})=\int_{0}^{z} \frac{ [v(\zeta)-(c+\frac{c^{*}}{\Gamma(a)})\zeta^{a-1}]}
{\left(z-\zeta \right) ^{a}} d\zeta =0
\end{equation*}
for all $z\in \mathbb{D}_{R}$ and for all $c,c^{*}\in \mathbb{C}.$ This implies that $v(z)=cz^{a-1}.$ Hence, $u(z)=cz^{a-1}$ are unique solutions of (3.1).  

\medskip\noindent\textbf{Lemma 3.2.} Under the conditions of Lemma 3.1, the following assertions are provided. 

\smallskip (i) $D_{z}^{a}I_{z}^{a}u(z)=u(z)$ for all $z\in \mathbb{D}_{R}.$  

\smallskip (ii) If $D_{z}^{a}u$ is continuous and integrable on $\mathbb{D}_{R},$ and $z^{a}D_{z}^{a}u$ is continuous and integrable on $\mathbb{U}_{R},$ then the equality
\begin{equation*}
I_{z}^{a}D_{z}^{a}u(z)=u(z)+cz^{a-1}     \ \  (c\in \mathbb{C})
\end{equation*}
holds for all $z\in \mathbb{D}_{R}.$   

\smallskip\noindent \textbf{Proof.}   (i) Let $u$ be continuous and integrable on $\mathbb{D}_{R}$. By using that $D_{z}^{a}=DI_{z}^{1-a}$ and
semi group property of the fractional integral we get 
\begin{equation}
D_{z}^{a}I_{z}^{a}u(z)=DI_{z}^{1-a}I_{z}^{a}u(z)=D_{z}I_{z}u(z)=u(z)  \tag{3.2}
\end{equation}
for all $z\in \mathbb{D}_{R}.$ The last equality in (3.2) holds for any continuous function $u$, 
since the integral $I_{z}$ is over the line segment. 

\smallskip (ii) At first, set 
\begin{equation}
v\left( z\right):=I_{z}^{a}\left(D_{z}^{a}u(z)\right) 
\tag{3.3}
\end{equation}
for all $z\in \mathbb{D}_{R}.$ Let us show that $u(z)=v(z)+cz^{a-1}$. If $D_{z}^{a}$ is applied to both sides of (3.3) and, after that, if the equality in (i) is used for the right side of the obtained equality, then the equality
\begin{equation*}
D_{z}^{a}v(z)=D_{z}^{a}u(z).
\end{equation*}
is obtained. 
Hence, if the linearity of the fractional derivative is considered, then it is derived that  
\begin{equation*}
D_{z}^{a}\left(\left(u-v\right)(z)\right) =0.
\end{equation*}
From the assumption, it follows that $u-v$ is continuous and integrable on $\mathbb{D}_{R}.$  Since the condition of Lemma 3.1 is satisfied, from the last equation it is clear that $u(z)=v(z)+cz^{a-1}.$ This is the desired result.  \\

\noindent\textbf{Remark 3.3.} (i) The unique continuous solution on $\mathbb{U}_{R}$ among solutions of the equation (3.1) is $u=0.$ Therefore, if one takes $u\in \mathbb{C}^{0}\left(\mathbb{U}_{R}\right)$ in hypotheses in Lemma 3.1 and Lemma 3.2, then $c$ in Lemma 3.2 (ii) has to be equal to zero. So, it is clear that $I_{z}^{a}$ is inverse of the operator $D_{z}^{a}$ with the domain $\mathbb{C}^{0}(\mathbb{U}_{R}).$  

\smallskip(ii) The equalities in Lemma 3.2 are also valid on $\overline{\mathbb{U}}_{R}$ and $\mathbb{D}^{*}_{R}=\big\{ z\in\overline{\mathbb{U}}_{R}: -\pi<\arg z\leq\pi \big\}.$  \\

 Now, we suppose that the condition (I) is satisfied and that the problem (1.1) has a solution $u\in \mathcal{B}.$ Then, $f(z,u(z))$ is continuous and integrable on $\mathbb{D}^{*}.$ By considering this fact in (1.1), it can be seen that $D_{z}^{a}u(z)$ is continuous and integrable on $\mathbb{D}.$ Hence, if $I_{z}^{a}$ is applied to the both sides of the equation (1.1) and, after that if Lemma 3.2 (ii) and Remark 3.3 (i) are used in the obtained equation, then the integral equation  
\begin{equation*}
u\left( z\right) \ = \frac{1}{\Gamma \left( a\right) }\  \int_{0}^{z}%
\frac{\ f\left( \zeta ,u\left( \zeta \right) \right) }{\left( z-\zeta
\right) ^{1-a}}d\zeta \eqno{(3.4)}
\end{equation*}
holds for all $z\in \overline{\mathbb{U}}.$ In addition to this, the equality (3.4) holds for $z=0,$ 
provided that the condition (II) is satisfied. \\

Consequently, the following Lemma can be deduced from the above explanations.  

\noindent\textbf{Lemma 3.4.} Let the conditions (I)-(II) be satisfied. If $u\in\mathcal{B}$ with $u(0)=b$ , then $u$ is a solution of the problem (1.1) if, and only if, $u$ satisfies the Volterra-type integral equation in (3.4).  \\

\noindent\textbf{Remark 3.5.} If the condition (II) does not hold, the contradiction can be obtained as follows:  
$$b=u(0)=\lim_{z\rightarrow 0}u(z)=\frac{1}{\Gamma(a)}\lim_{z\rightarrow 0} \int_{0}^{1}\frac{(zt)^{a}\ f\left(zt,u\left(zt\right) \right) }{t^{a}\left( 1-t\right)^{1-a}}dt\neq b.$$   

In the following, we first prove the existence and uniqueness of the solution for (1.1) with the initial condition $u(0)=0,$ and then for that subjected to initial condition $u(0)=b.$

\smallskip\noindent\textbf{Theorem 3.6.}    
Let the condition (I) be satisfied. Moreover, we assume that there exist a fixed natural number $n_{0}\geq 1,$ a non-negative real number $c$ and a function $g\in\mathcal{B}^{0}$ such that the following inequality holds for all $(z,t)\in\overline{\mathbb{U}}\times \mathbb{C}:$
$$\left \vert z^{a}f\left( z,t\right) \right \vert \leq c\left|t\right|^{n_{0}}+\left|g(z)\right|. \eqno{(3.5)}$$

\noindent Then there exists a $R\in (0,1]$ such that the problem (1.1) with $b=0$ has at least one solution $u\in \mathcal{B}_{R}^{0}.$ \\
 
\smallskip\noindent\textbf{Proof.} Suppose that $u\in \mathcal{B}^{0}.$  From our assumptions it follows that the hypotheses of Lemma 3.4 are satisfied. Thus, the problem (1.1) is equivalent to the integral equation (3.4). If the operator $P$ is defined as   
\begin{equation*}
Pu\left(z\right)=\frac{1}{\Gamma\left(a\right) }\int_{0}^{z}%
\frac{f\left(\zeta ,u\left( \zeta \right) \right) }{\left(z-\zeta
\right)^{1-a}}d\zeta, 
\end{equation*}
then $P$ is an operator from $\mathcal{B}^{0}$ to $\mathcal{B}^{0}.$  Hence, the fixed points of $P$ in $\mathcal{B}^{0}$ coincide the solutions of the problem (1.1). Thus, it is sufficient to prove the existence of the fixed points of the operator $P.$ For the proof, it is shown that the all conditions of Theorem 2.3 are satisfied. \\

Consider first the bounded, closed and convex subset $B_{r}$ of $\mathcal{B}_{R}^{0}$ given by $B_{r}=\big \{u\in\mathcal{B}_{R}^{0}:\left \Vert u\right \Vert_{\mathcal{B}_{R}^{0}}\leq r\big\}$ with the fixed $r>0$ and $R\in(0,1],$ and let us see that there exist a suitable $R$ such that 
\begin{equation*}
P(B_{r})\subseteq B_{r} \tag{3.6}
\end{equation*}
is satisfied.  \\

Now, since for all $(z,t)\in\overline{\mathbb{U}}\times \mathbb{C}$ the inequality (3.5) is satisfied
, it follows that 
$$\left \vert z^{a}f\left( z,u(z)\right) \right \vert \leq c\left|u(z)\right|^{n_{0}}+\left|g(z)\right|$$
for all $z\in\overline{\mathbb{U}}$ and for all $u\in\mathcal{B}^{0}.$ \\ 

By using Schwarz Lemma in the above inequality one can see that 
$$\left \vert z^{a}f\left( z,u(z)\right) \right \vert \leq cr^{n_{0}}\left|z\right|^{n_{0}}+M_{g}\left|z\right| $$ 
is satisfied  on $\overline{\mathbb{U}}$ for all $u\in\mathcal{B}^{0}$  with $\left|u(z)\right|\leq r$ and for $g\in\mathcal{B}^{0}$ with $\left|g(z)\right|\leq M_{g}.$ \\

From above inequality  it is derived that 
\begin{align*}
\left \vert Pu\left( z\right) \right \vert &\leq \frac{1}{\Gamma \left(a\right)} 
\sup_{z\in \overline{\mathbb{U}}_{R}}\left[\int_{0}^{1}\frac{cr^{n_{0}}\left|z\xi\right|^{n_{0}}}{\xi^{a}\left(1-\xi%
\right)^{1-a}}d\xi+M_{g}\int_{0}^{1}\frac{\left|z\xi\right| }{
\xi^{a}\left(1-\xi \right)^{1-a}}d\xi \right] \\
&\leq cr^{n_{0}}R^{n_{0}}\frac{\Gamma \left(n_{0}+1-a\right)}{\Gamma \left(n_{0}+1\right)}+M_{g}R\Gamma \left(2-a\right)
\end{align*}
for all $|z|\leq R\leq 1$. Therefore, it can be easily seen that (3.6) holds for a suitable $0<R=R(r,n_{0},M_{g},a)<1.$ \\

It remains to show that $P$ is a continuous operator on  $B_{r}$ and $P(B_{r})$ is an equicontinuous set of $\mathcal{B}_{R}^{0}.$ For the continuity of $P$ on $B_{r},$ it is supposed that $\left \{u_{n}\right \}_{n=1}^{\infty}\subset B_{r}$ is a sequence with $u_{n}\stackrel{\mathcal{B}_{R}^{0}}{\rightarrow} u$ as $n\rightarrow \infty.$ Then, it is clear that $u_{n}$
converges uniformly to $u\in B_{r},$ since $B_{r}$ is a  closed subset of $\mathcal{B}_{R}^{0}.$
By using the uniform continuity of $z^{a}f(z,t)$ on $\overline{\mathbb{U}}_{R}\times \overline{\mathbb{U}}_{r}$ ($\overline{\mathbb{U}}
_{r}:=\left \{ \nu \in \mathbb{C}:\left|\nu \right|\leq r \right \}$) and uniform
convergence of ${u_{n}}$ to the function $u$ on $\overline{\mathbb{U}}_{R},$ one can conclude that
\begin{align*}
&\left \|Pu_{n} -Pu\right \|_{ \mathcal{B}^{0}_{R_{0}}}\\
&\ \ =\sup_{z\in \overline{\mathbb{U}}_{R_{0}}} \left \vert \frac{1}{\Gamma \left(
a\right) }\  \int_{0}^{z}\frac{ \left[f\left( \zeta ,u_{n}(\zeta
)\right) -f\left( \zeta ,u(\zeta )\right)\right] }{\left( z-\zeta \right)
^{1-a}}d\zeta \right \vert \\
&\ \ \leq \frac{1}{\Gamma \left( a\right)} \sup_{\xi z\in \overline{%
\mathbb{U}}_{R_{0}}} \int_{0}^{1}\frac{ \left|(\xi z)^{a}f\left(\xi z,u_{n}(\xi z)\right) -(\xi z)^{a}f\left( \xi z ,u(\xi z)\right)\right|}{\xi^{a}\left( 1-\xi \right) ^{1-a}} d\xi \rightarrow 0
\end{align*}%
as $n\rightarrow \infty.$  

Now, let us show that $P(B_{r})$ is an equicontinuous set of $\mathcal{B}_{R}^{0}.$ 
Since all $u\in B_{r}$ are uniformly continuous on $\overline{\mathbb{U}}_{R}$ and $z^{a}f(z,t)$ is
uniformly continuous on $\overline{\mathbb{U}}_{R}\times \overline{\mathbb{U}}_{r},$ 
then $z^{a}f(z,u(z))$ is also uniformly continuous on $\overline{\mathbb{U}}_{R}.$ Therefore, for given $\epsilon>0$ there exists a $%
\delta=\delta(\epsilon)>0$ such that 
$$\left \vert z_{1}^{a}f(z_{1},u(z_{1}))-z_{2}^{a}f(z_{2},u(z_{2}))\right \vert <\frac{\epsilon}{\Gamma(1-a)},$$
for all $z_{1},z_{2}\in \overline{\mathbb{U}}_{R}$ satisfying $\left|z_{1}-z_{2}\right|<\delta.$ From here, one can conclude that 
\begin{align*}
&\big|Pu\left( z_{1}\right) -Pu\left( z_{2}\right) \big| \\
&\ \ \ \leq \frac{1}{\Gamma \left( a\right)}\int_{0}^{1} \frac{\left \vert
(\xi z_{1})^{a}f\left( \xi z_{1},u\left( \xi z_{1}\right) \right) -(\xi
z_{2})^{a}f\left( \xi z_{2},u\left( \xi z_{2}\right) \right) \right
\vert}{\xi^{a} \left( 1-\xi \right) ^{1-a}} d\xi \  \\
&\ \ \ <\Gamma(1-a)\frac{\epsilon}{\Gamma(1-a)}=\epsilon,
\end{align*}
since $\left|\xi z_{1}-\xi z_{2}\right|<\delta.$ So, it is obtained that $P(B_{r})
$ is an equicontinuous set of $\mathcal{B}_{R}^{0}.$  

\smallskip Consequently, as a consequence of Theorem 2.3, one can say that the operator $P$ has at least one fixed point in $\mathcal{B}_{R}^{0}$ for a $R\in (0,1]$ given above, and it is also a solution of the problem (1.1).  \\

\noindent\textbf{Remark 3.7.} Schauder's fixed point theorem is applicable to prove the existence of local continuous solution for the problems with FDE in many researches. However, the existence of local desired solution of the problem (1.1)  can be proved only for a subclass of functions $f$ satisfying the conditions (I) and (II). Indeed, from the condition (I) one can suppose that 
\begin{align*}
\left|z^{a}f(z,t) \right|\leq M \ \ (\forall (z,t)\in\overline{\mathbb{U}}\times\overline{\mathbb{U}}_{r}), 
\end{align*}
and by using this inequality one can write 
\begin{align*} 
\sup_{z\in \overline{\mathbb{U}}_{R}}\left \vert Pu\left( z\right) \right \vert\leq \frac{M}{\Gamma \left(a\right)}
\int_{0}^{\left|z\right|}\frac{1}{\left|\zeta\right|^{a}\left|z-\zeta\right|^{1-a}}\left|d\zeta\right|\leq M \Gamma(1-a)\tag{3.7} 
\end{align*} 
for all $u\in B_{r}$ and for all $z\in\overline{\mathbb{U}}_{R}$ with an arbitrary $R\in \left(0,1\right].$ Hence, it must be $M \Gamma(1-a)\leq r$ in order that the condition $P(B_{r})\subseteq B_{r}$ is satisfied. This indicates that the function $f(z,t)$ has to satisfy the following inequality:   
\begin{align*}
\left|z^{a}f(z,r)\right|\leq \frac{r}{\Gamma(1-a)} \ \ \  (\forall z\in\overline{\mathbb{U}}), 
\end{align*} 
which means that $f$ increases not faster than a linear function of $r.$ Moreover, this inequality is a particular case of the inequality in Theorem 3.6, when $c:=\frac{1}{\Gamma(1-a)},$ $n_{0}=1$ and $g(z)\equiv 0.$ In spite of this, the solution of the considered problem exists on whole $\mathbb{U},$ since the inequality (3.7) holds for all $z\in\overline{\mathbb{U}}_{R}$ with an arbitrary $R\in \left(0,1\right].$ \\

 Theorem 3.6 does not imply that the problem (1.1)  with $b=0$ admits a unique solution in $\mathcal{B}_{R}^{0}.$  Indeed, if $f(z,u):=\frac{z^{-a}}{\Gamma \left(2-a\right)}u$ with the fixed $a\in (0,1)$ in this problem, then it admits the solutions $u(z)=c_{*}z,$ where $c_{*}\in\mathbb{C}.$ Now, we give in the following theorem which implies not only existence but also uniqueness of the desired solution for the considered problem.

\smallskip\noindent\textbf{Theorem 3.8.} Let the conditions (I) and (II) with $b=0$ be satisfied.  
Moreover, assume that there exists a
constant $\kappa< 1/\Gamma \left(2-a\right)$ such that 
\begin{equation}
\left \vert f\left( z,\eta \right) -f\left( z,\nu \right) \right \vert <%
\frac{\kappa}{\left|z\right|^{a}}\left \vert \eta -\nu \right \vert 
\tag{3.8}
\end{equation}%
for all $z\in \mathbb{D^{*}} $ and for all $\eta $,$\nu \in \mathbb{C}.$
Then the problem (1.1)  with $b=0$ has a unique solution in $\mathcal{B}_{0}.$  

\smallskip\noindent\textbf{Proof.}  It is supposed that $u\in\mathcal{B}_{0}.$  Since the conditions of Lemma 3.4 are satisfied, 
the operator $P$ defined in proof of Theorem 3.6 can be considered. The operator $P$ is well defined. Indeed, if the inequality in
(3.8) is taken into account, then it is obtained that the inequality 
$$\left|z^{a}f(z,\eta)\right|<\kappa \left|\eta\right|+\left|z^{a}f(z,0)\right|$$ 
for all $z\in \overline{\mathbb{U}}$ and for all $\eta \in \mathbb{C}.$ 
By using this inequality 
\begin{align*}
\left \vert Pu\left( z\right) \right \vert\leq \Gamma \left( 1-a\right)\left(\kappa \left \|u\right \|_{\mathcal{B}_{0}}+\sup_{z\in \overline{\mathbb{U}}}\left \vert
z^{a}f(z,0)\right \vert \right)
\end{align*}
can be obtained.  

\smallskip Hence, the fixed points of the operator $P$ coincide with the solutions of the considered problem. Then, it is sufficient to show that this \linebreak operator has a unique fixed point in $\mathcal{B}_{0}.$  To do this by virtue of Banach fixed point theorem,
it is enough to see that the operator $P$ is a contraction.  

\smallskip Now, let us suppose that $u$ and $u_{0}$ are the arbitrary 
elements of $\mathcal{B}_{0}.$ Then, $u-u_{0}\in \mathcal{B}_{0}.$ By using Schwarz's Lemma the following chain of inequalities can be obtained: 
\begin{align*}
\left \vert Pu\left(z\right)-Pu_{0}\left(z\right) \right \vert &\leq 
\frac{1}{\Gamma \left( a\right)}\int_{0}^{\left|z\right|}\frac{\left
\vert f\left(\zeta ,u(\zeta )\right) -f\left( \zeta ,u_{0}(\zeta
)\right)\right| }{ \left|z-\zeta \right| ^{1-a}}\left|d\zeta \right| \\
&\leq\frac{\kappa}{\Gamma \left( a\right) }\  \int_{0}^{\left|z\right|}
\frac{\left|\zeta \right| \left|\frac{\left(u -u_{0}\right)(\zeta)}{\zeta}
\right| }{\left|\zeta \right|^{a}\left|z-\zeta \right| ^{1-a}}
\left|d\zeta \right| \\
&\leq \kappa \Gamma \left( 2-a\right)\left \Vert u-u_{0}\right \Vert _{\mathcal{B}_{0}}.
\end{align*}
Since $\kappa<1/\Gamma \left( 2-a\right),$ the above inequality implies that $P$ is a contraction operator. As a consequence of Banach fixed point theorem, one can say that there exists
a unique fixed point of the operator $P$ in the space $\mathcal{B}_{0}.$ Consequently, the considered problem has a unique solution $u$ in $\mathcal{B}_{0}.$ \\

\noindent\textbf{Remark 3.9.} (i) If the technique related to Schwarz's Lemma wasn't used in the proof of Theorem 3.8, then it would be obtained the inequality: 
\begin{align*}
\left \vert Pu\left( z\right)-Pu_{0}\left( z\right) \right \vert<\kappa_{1} \Gamma \left( 1-a\right)\left \Vert u-u_{0}\right \Vert _{\mathcal{B}_{0}}
\end{align*} 
for any $u,u_{0}\in\mathcal{B}_{0}.$ For $P$ to be contraction operator we require that $\kappa_{1}<1/\Gamma \left( 1-a\right).$ Since
$\kappa_{1}<\kappa_{2},$ it is clear that using Schwarz Lemma in the proof of the related theorem contributed the function $f$ to be in a larger class of functions. \\  

(ii)  Under the conditions in Theorem 3.6 or Theorem 3.8, the solutions of the considered problem need not to be univalent in related domain, unless some more conditions on the function $f$ are imposed. Indeed, if $f(z,t):=cz^{-a}t$  with  $c=\frac{\Gamma(n+1)}{\Gamma(n+1-a)}$ for fixed $n\in \mathbb{N}-\left\{1\right\}$ in this problem, then it has the non-univalent solutions $u(z)=c^{*}z^{n}$ for all $c^{*}\in\mathbb{C}.$ \\
 
The condition (II) is necessary for the equivalence of (1.1)  with $b=0$ and the Volterra type equation in (3.4). In the following it is shown that this condition is indispensable by proving that, without this condition, the considered problem has no analytic solution.

\smallskip\noindent\textbf{Proposition 3.10.} The condition $z^{a}f\left( z,0\right)|_{z=0}=0$ is necessary for the existence of analytic solution to the problem (1.1)  with $b=0.$
 
\smallskip\noindent\textbf{Proof.}  This proposition is proved by showing that there exist a contradiction. For this, let $a$ be fixed in $
(0,1)$ and let $f(z,u):=cz^{-a}u+dz^{-a}$ with any $d\in\mathbb{C}-\left\{0\right\},$ $c\in \mathbb{C}$ in the considered problem.  It is supposed that this problem admits an analytic solution $u$ on $\mathbb{U}.$ If one considers that $u$ is represented by a power-series expansion such as $u(z)=\sum^{\infty}_{k=1}a_{k}z^{k}$ on $\mathbb{U},$ then the following equality can be obtained: 
\begin{align*}
\sum^{\infty}_{k=1}a_{k}\frac{\Gamma(k+1)
}{\Gamma(k+1-a)}z^{k-a}=\sum^{\infty}_{k=1}ca_{k}z^{k-a}+dz^{-a}   \ \ (z\in\mathbb{D}).
\end{align*}
However, it is a contradiction since the above equality does not hold for $d\in\mathbb{C}-\left\{0\right\}$. Therefore, there is no analytic solution of this problem unless the condition $z^{a}f\left( z,0\right)|_{z=0}=0$ is satisfied.  \\

The existence and uniqueness results given above can be obtained for the problem (1.1)  with non-homogenous initial data. However, it is sufficient to focus on some differences in the following. 

\smallskip
\noindent\textbf{Remark 3.11.} Now, we consider the problem (1.1)  with non-homogeneous initial data. 
By replacing $u(z)$ by $v(z)+b$ in the problem (1.1) , the following problem is obtained: 
\begin{equation*}
\begin{array}{rcl}
D_{z}^{a}v(z)&=&h(z,v(z)+b)  \\ 
v(0)& = & 0.
\end{array}
\end{equation*}
where $h(z,v(z)+b)=-\frac{b}{\Gamma(1-a)}z^{-a}+f(z,v(z)+b).$ If it is supposed that the condition (I) and the condition 
\begin{align*}
\left \vert z^{a}h\left(z,t\right) \right \vert \leq c\left|t-b\right|^{n_{0}}+\left|g(z)\right|, \ \ (n_{0}\in \mathbb{N}, c\geq 0, g\in\mathcal{B}^{0})  \tag{3.9}
\end{align*}
instead of the inequality (3.5) in Theorem 3.6 are satisfied, then it can be shown, under these conditions, that 
the problem just above has at least one solution $v\in\mathcal{B}^{0}_{R}$ for a suitable $0<R\leq 1.$ So, from the equivalence of these problems, one can say that the problem (1.1)  admits at least one solution $u\in\mathcal{B}_{R}$ with $u(0)=b.$
\smallskip

Furthermore, if one changes the condition $z^{a}f(z,0)|_{z=0}=0$ in Theorem 3.8  with $z^{a}f(z,b)|_{z=0}=\frac{b}{\Gamma(1-a)}$ by keeping the other conditions of this theorem same, then one can prove, under these new conditions, that the problem (1.1)  admits unique solution $u\in\mathcal{B}$ with $u(0)=b.$ \\

Now let us consider the problem 1.2. Using Lemma 3.2, it can be easily seen that every solution $u\in \mathcal{B}$ with $u(0)=b$ of the 
problem (1.2)  is  also a solution of the following fractional integral equation:
$$u\left( z\right)=b+ \frac{1}{\Gamma \left( a\right) }\  \int_{0}^{z}\frac{\ f\left( \zeta ,u\left( \zeta \right) \right) }{\left( z-\zeta
\right) ^{1-a}}d\zeta $$
and vice versa. 

\smallskip
In the view of the results in Remark 3.11 and by using the same ways in proofs of Theorem 3.6 and Theorem 3.8, we give the following existence and uniqueness results without proofs:  

\smallskip

\noindent\textbf{Theorem 3.12.}   Let the condition (I) be satisfied. In addition to this, we suppose that there exist a fixed natural number $n_{0}\geq 1,$ a non-negative real number $c$ and a function $h\in\mathcal{B}^{0}$ such that the following inequality holds for all $(z,t)\in\overline{\mathbb{U}}\times \mathbb{C}:$
$$\left \vert z^{a}f\left( z,t\right) \right \vert \leq c\left|t-b\right|^{n_{0}}+\left|h(z)\right|. \eqno{(3.10)}$$

\noindent Then there exists a $R\in (0,1]$ such that the problem (1.2)  has at least one solution $u\in \mathcal{B}_{R}$ with $u(0)=b.$\\

\smallskip

\noindent\textbf{Theorem 3.13.} Let the condition (I) and $z^{a}f(z,b)\big|_{z=0}=0$ be satisfied.  
Furthermore, assume that there exists a constant $\kappa< 1/\Gamma \left(2-a\right)$ such that 
$$\left \vert f\left( z,\eta \right) -f\left( z,\nu \right) \right \vert <%
\frac{\kappa}{\left|z\right|^{a}}\left \vert \eta -\nu \right \vert $$
for all $z\in \mathbb{D^{*}} $ and for all $\eta $,$\nu \in \mathbb{C}.$
Then the problem (1.2)  admits a unique solution $u\in \mathcal{B}_{R}$ with $u(0)=b.$

\medskip\noindent{\textbf{\large 4. Existence and Uniqueness Results for Problem (1.3)}} 

\smallskip

\noindent In this section, we mention about the existence and uniqueness of the solution to the problem (1.3). For this, in the rest of this section, we suppose that it is given a function $f(x,y)$ defined in $[0,1]\times\mathbb{R}$ and can be analytically continued to a function $f(z,t)$ defined in $\overline{\mathbb{U}}\times\mathbb{C}.$ Then, as the results of theorems given above, we obtain the following theorems for the problem (1.3)  with Riemann-Liouville (R-L) or Caputo derivative:    

\smallskip\noindent\textbf{Theorem 3.14.}  If the function $f(z,t)$  which is the analytic continuation of  given function $f(x,y)$  satisfies the conditions (I) and the inequality (3.9), and if 
$$f(x,y)=\Re\big(f(z,t)\big) \eqno{(3.11)}$$ 
holds, then the problem (1.3)  with R-L derivative  has at least one continuous solution on the interval $[0,R_0]$ with a suitable $R_{0}\leq 1.$ 

\smallskip\noindent\textbf{Proof.}  Suppose that the function $f(z,t)$ satisfies the conditions (I)-(II), which is the analytic continuation of  given function $f(x,y)$ right hand side of (1.3). As a consequence of Remark 3.11, one can say that there exists a solution $u\in \mathcal{B}_{R}$ with $u(0)=b,$ i.e. $u$ satisfies the equation (1.1) on $\overline{\mathbb{U}}_{R}$ and therefore on $[0,R].$ In that case, if we take the real part of the both sides of equation (1.1) and, if we use the assumption and the equality (2.3), then we have 
$$\mathcal{D}^{a}\left(\Re \left\{u(z)\right\}\right)= f\big(x,\Re \left\{u(z)\right\}\big) ,$$
which means that $\Re \left\{u(z)\right\}$ is a continuous solution of (1.3). \\

\noindent\textbf{Theorem 3.15.}  Assume that the function $f(z,t)$ which is the analytic continuation of given function $f(x,y)$ fulfills the conditions (I)-(II) and the inequality (3.8). Moreover, if the equality (3.11) is satisfied, then the problem (1.3) with R-L derivative has a unique continuous solution on the interval $[0,1],$ which is real analytic on $(0,1).$ \\

\noindent\textbf{Corollary 3.16.}  Suppose that the function $f(z,t)$ which is the analytic continuation of given function $f(x,y)$  satisfy the conditions (I)-(II) and, let $z^{a}f(z,t)$ be a linear function in $z$ and $t$ . Then the problem (1.3) with R-L derivative admits a unique solution continuous on the interval $[0,1]$ and real analytic on $(0,1).$ \\

As the consequences of Theorems 3.12-3.13 and the well-known relation $^{C}D^{a}(u(x))=D^{a}(u(x)-u(0))$
for the differentiable function $u$ between Caputo and Riemann-Liouville derivatives, one can conclude the following theorems and corollary by the same way used for obtaining Theorems 3.14-3.15 and Corollary 3.16.

\smallskip\noindent\textbf{Theorem 3.17.}  Let the function $f(z,t)$ which is the analytic continuation of given function $f(x,y)$ satisfy the conditions (I) and, let the inequality (3.10) be satisfied for a fixed natural number $n_{0}\geq 1,$ a non-negative real number $c$ and a function $h\in\mathcal{B}^{0}.$ Furthermore, if the equality (3.11) 
holds, then the problem (1.3)  with Caputo derivative possesses at least one solution continuous on the interval $[0,R_0]$ and real analytic on $(0,R_0)$ for a suitable $R_{0}\leq 1.$ \\

\noindent\textbf{Theorem 3.18.}  Assume that the function $f(z,t)$ which is the analytic continuation of given function $f(x,y)$ fulfills the condition (I),  $z^{a}f(z,b)\big|_{z=0}=0$ and the inequality (3.8). Furthermore, if the equality (3.11) 
is fulfilled, then the problem (1.3)  with Caputo derivative admits a unique continuous solution on the interval $[0,1],$ which is real analytic on $(0,1).$ \\

\noindent\textbf{Corollary 3.19.}  Suppose that the function $f(z,t)$ which is the analytic continuation of given function $f(x,y)$ satisfies the conditions (I) and $z^{a}f(z,b)\big|_{z=0}$ $=0$. Moreover, let $z^{a}f(z,t)$ be a linear function in $z$ and $t.$ Then the problem (1.3)  with Caputo derivative has a unique solution continuous on the interval $[0,1]$ and real analytic on $(0,1).$

\end{document}